\documentclass{article}

\usepackage{ucs}
\usepackage{pdfpages}
\usepackage[T1]{fontenc}
\usepackage{amssymb}
\usepackage{amsmath}
\usepackage{amsthm}
\usepackage{verbatim}
\usepackage{multirow}
\usepackage{tabularx}
\usepackage[linesnumbered,ruled]{algorithm2e}

\usepackage{here}
\usepackage{hyperref}
\usepackage{cite}

\usepackage[utf8]{inputenc}
\usepackage[english]{babel}
\usepackage{graphicx}

\usepackage{microtype}

\newtheorem{theorem}{Theorem}
\newtheorem{lemma}{Lemma}
\newtheorem{corollary}{Corollary}
\newtheorem{definition}{Definition}
\newtheorem{remark}{Remark}
\newtheorem{example}{Example}
\newtheorem{proposition}{Proposition}

\usepackage[all]{xy}

\usepackage{ifthen}
\newboolean{showproofs}
\setboolean{showproofs}{true}

\begin{document}

\newpage

\title{\bf How to avoid collisions in 3D-realizations for moving graphs}

\author{Jiayue Qi \thanks{Johannes Kepler University Linz, Doctoral Program 
"Computational Mathematics" (W1214).} 
\thanks{Johannes Kepler University Linz, Research Institute for Symbolic Computation.}}

\date{}
\vspace{0.2cm}

\maketitle
{\centering\footnotesize\bf\em To my dearest Grandma, Huiqin Dong.\par}

\begin{abstract}

If we parameterize the positions of all vertices of a given graph in the plane such that
distances between adjacent vertices are fixed, we obtain a moving graph. 
An L-linkage is a realization
of a moving graph in 3D-space, by representing edges using horizontal bars and 
vertices by vertical sticks. Vertical sticks are parallel revolute joints, 
while horizontal bars are links connecting them. 
We give a sufficient condition for a moving graph to have a collision-free
L-linkage. Furthermore, we provide an algorithm guiding the construction of such a 
linkage when           
the moving graph fulfills the sufficient condition, via computing a height function for the edges 
(horizontal bars).
In particular, we prove that any Dixon-1 moving graph has a collision-free L-linkage and 
no Dixon-2 moving graphs have collision-free L-linkages, where Dixon-1 and Dixon-2 moving 
graphs
are two classic families of moving graphs. 

\end{abstract}

\section{Introduction}\label{sec:introduction}
In this section, we introduce the problem background and some related existing work.
Given a graph, if we parameterize the position of each vertex in the plane so that 
any distance between 
adjacent vertices is constant,
then this graph is called a moving graph. Notice that some moving graphs are 
induced by rigid motions, hence they do not necessarily ``move'' in the general sense. 
Examples of moving graphs can be found in 
\cite{moving_graph, Dixon}.
We also provide more examples of moving graphs in the later sections.

\begin{figure}[H]
\centering
\includegraphics[width=0.4\linewidth]{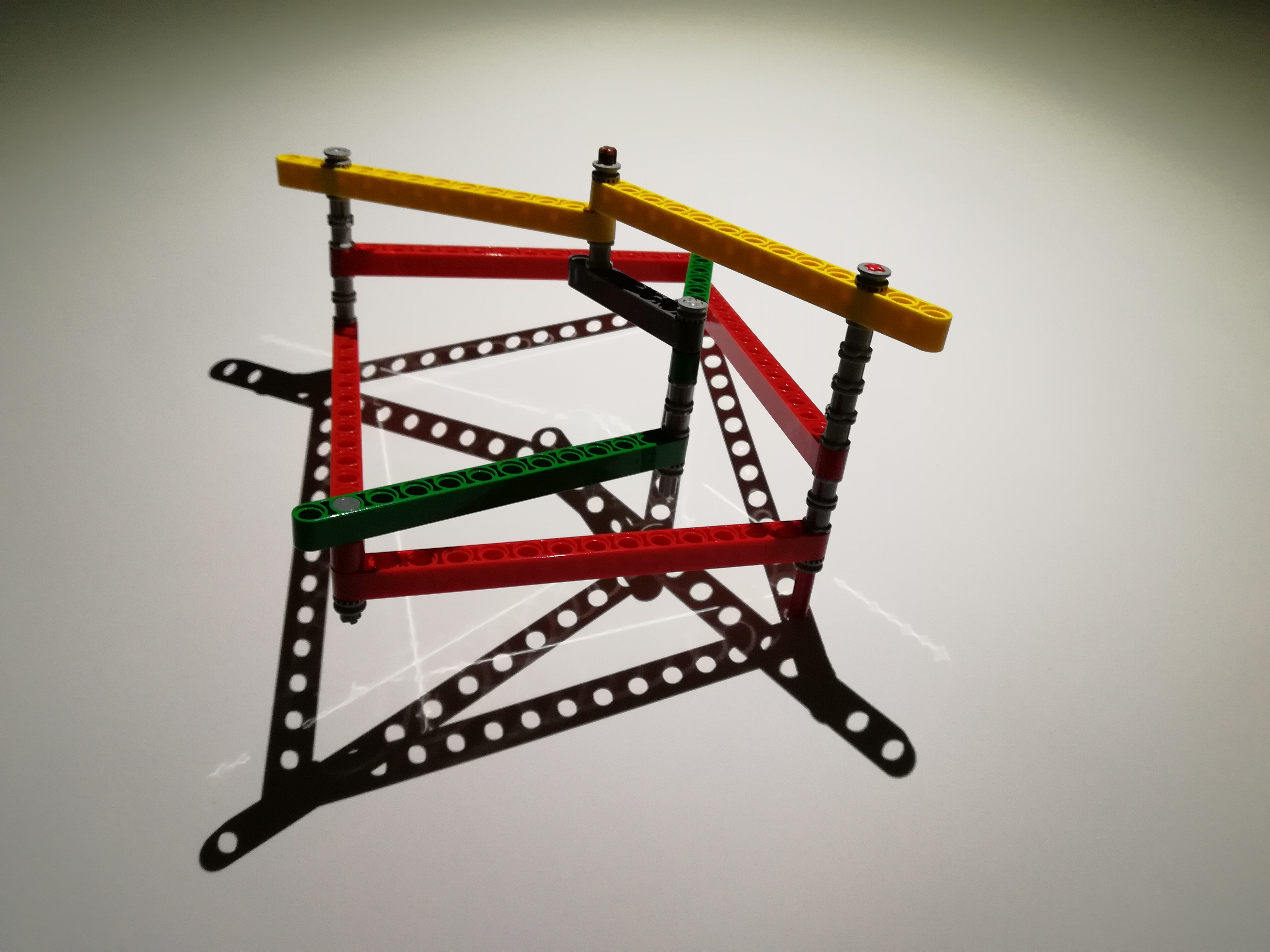}
\caption{This is the realization of an L-model (i.e. an L-linkage) of a moving graph with the underlying 
graph structure being $K_{3,3}$, where
$K_{3,3}$ denotes the complete bipartite graph with two independent sets and each independent set contains three vertices.
The edges correspond to horizontal bars with holes, while the vertices correspond to vertical sticks. For 
detailed explanation of complete bipartite graphs see Section \ref{sec:collision_detection}.}
\label{fig:photo}
\end{figure}

We realize a moving graph in 3D-space with some bars and sticks. Horizontal bars 
realize edges while vertical 
sticks realize vertices. We assign a different height to each horizontal bar. 
Horizontal bars are parallel and located at pairwise distinct heights and hence would not
touch each 
other. Vertical sticks connect them, by going through the holes at the ends of horizontal 
bars. We call this kind of linkages L-linkages, or realizations of L-models. 
The concept of L-model is just the mathematical description of such linkages, which we will 
define precisely later on, in mathematical language.
Figure \ref{fig:photo} shows an L-linkage.

 One observes a collision of an L-linkage at the time when a vertical stick (corresponds to
 vertices of the corresponding moving graph of the linkage) hits a horizontal bar  
 (corresponds to
 edges of the corresponding moving graph of the linkage).
 The collisions depend on the choice of heights of the 
horizontal bars. Note that edge crossing is possible in some moving graphs. 
However, we set up horizontal bars all in different heights, which naturally
avoids edge crossing in the corresponding L-linkage.
We give a sufficient condition for such a linkage to be collision-free with respect
to some height arrangement of edges.
And we provide an algorithm for constructing an L-linkage in a collision-free way
when the criterion is fulfilled. This condition can be helpful for the construction of 
the linkage --- checking whether it can accomplish a full motion.  
Dixon-1 moving graphs and Dixon-2 moving graphs
are two classic families of moving graphs. With this result we prove that there exists a 
collision-free L-linkage
for every Dixon-1 moving graph. In addition, we prove that Dixon-2 moving graphs cannot be realized 
as a collision-free L-linkage.

  Abel et al. showed that any polynomial curve can be traced by a non-crossing linkage 
  \cite{who_needs_crossings}. One of the differences between our result and theirs is
that we consider a motion (i.e. the traces of all vertices),
not just tracing a curve of some single vertex.
Another relevant result is by Gallet et al. \cite{prescribed_motion}. They provide an algorithm which produces a linkage tracing 
any given planar rational curve without collisions. Our work differs from them in the sense that
our algorithm deals with arbitrary
moving graphs, not just those related to their linkages; also, they do not take into consideration 
the moving graphs 
Dixon-1 and Dixon-2. The result in \cite{deployable} provides a method to better avoid collision 
for curved-bar-type linkages, while we focus on collision avoidance for straight-bar linkages.
Collision-free path planning of planar linkages is also discussed in \cite{path_planning}.
Our result can also be helpful in the collision-free path planning. Given a path, a corresponding 
linkage which draws that path can be designed with some method \cite{prescribed_motion}. Then, our algorithm can check whether the linkage
can be realized without collision, which contains the ability of checking whether the planned path can be 
realized in a collision-free way. The result in \cite{6R} also considers collision-avoidance. They
introduce a specific design of a 6R-linkage in such a way that it does not collide with itself or
other equipment attached to it. Compare with their result, we consider the collision-free design for another type of linkages, namely L-linkages.

\section{Problem statement}
In this section, we describe the problem on which this paper focuses.
For a clear problem statement we need to introduce some definitions first.
These definitions help us build up the mathematical model precisely, behind
the intuitional idea. If we view the linkage in Figure \ref{fig:photo} from 
above, we get some idea of a ``moving graph''.
\begin{definition}
Let $G=(V,E)$ be a connected graph. Let $F=\{f_v\}_{v\in V}$ be a set of continuous functions
$f_v:\mathbb{R}\rightarrow \mathbb{R}^2$ such that the function $\|f_u(\cdot)-f_v(\cdot)\|$
is a positive constant for every edge $\{u,v\}$ in $E$. The pair $(G,F)$ is called 
a {\bf moving graph}.
\end{definition}

We can fix the position of some vertex $v$ and put one of its adjacent vertices $v_1$ on the $x$-axis. 
Let $\lambda_{i,j}$ be the fixed edge length for edge $\{i,j\}$. We can consider the equation system which contains
equation $\|f_u(t)-f_v(t)\|=\lambda_{u,v}$ for each $\{u,v\}\in E$.
 If the corresponding solution set (in $\mathbb{R}$) has infinite cardinality, we say that the given moving graph is {\em mobile}.
 We also say that the parameterizations give a {\em mobile realization} of the underlying graph.
 The solution set is also called a configuration space. We can talk about the 
 dimension of this space. The space has a positive dimension if and only if the given moving graph is mobile. 
 When the solution set has finite cardinality, the given moving graph is {\em rigid}.
 For our collision problem, it only makes sense to consider the mobile moving graphs.

A moving graph $M=(G,F)$ is {\itshape finite} if $G$ is a finite graph. 
The background problem that we are interested in is how to 
avoid collisions for L-linkages, hence we focus only on finite moving graphs in this paper.

We define a collision in a moving graph as when a vertex ``hits'' an edge.
\begin{definition}\label{def:collision_pairs}
 Let $M=(G,F)$ be a moving graph, where $G=(V,E)$. Vertex {\bf $w\in V$ collides 
 with edge $\{u,v\}\in E$} if and only if $w\notin \{u,v\}$ and
 there exists $t\in \mathbb{R}$ such that $f_w(t)$ lies on the line segment defined by points $f_u(t)$ and $f_v(t)$.
 Then $(w,\{u,v\})$ is called a {\bf collision pair} of $M$. 
 We denote by $CP_M$ the set of all collision pairs 
 of~$M$. 
\end{definition}

Note that by definition, a vertex never collides with any edge containing this vertex. 
If we assign pairwise distinct (integer) height values for the edges, we obtain an L-model. 
This is the mathematical language for the Lego linkages (L-linkages) described in Section~\ref{sec:introduction}.

\begin{definition}
Let $M=(G,F)$ be a moving graph, where $G=(V,E)$. An {\bf L-model} of $M$ is a pair $(M,h)$, where $h:E\rightarrow \mathbb{Z}$ is an injective 
function assigning to each edge of $G$ an integer height value. We call $h$ the {\bf height function} of $M$.
\end{definition}
An {\em L-linkage} is a 3D-realization of an L-model, as described in Section \ref{sec:introduction}.
Although our original ambition is to avoid collisions for L-linkages, L-model is 
a concept in mathematical set-up that can help us better, in developing the whole theory. 
Hence, we focus more on this concept throughout the paper. 

In some moving graphs, two (or more) edges can stay intersected for any $t\in \mathbb{R}$. 
However, we can prevent such situation being transferred to a collision 
problem in the corresponding L-linkage, simply by placing all edges in pairwise distinct heights.
Because we define the height function to be injective, so this would not be a problem in our setting.
Therefore, we do not specifically consider this situation in the remaining context.

When we are given
an L-model, we can realize the edges by horizontal bars by fixing the height of one edge and then placing the others accordingly. 
And we connect the end of those edges that are incident at the same vertex by a vertical stick going through the holes on the horizontal 
bars. In this way, we obtain the corresponding L-linkage of the given L-model. Note that during this process, only the relative height
values matter. We can obtain a same L-linkage from two different L-models, as long as the relative height between any pair of 
edges coincides in these two L-models.
On the contrary, we can also
find out the height function of one corresponding L-model of a given L-linkage, by pre-fixing the height of an edge in the linkage.

Collision pairs are the ``collisions'' in a moving graph; also, they are the only reason for the corresponding L-linkage to have collisions.
When we realize the moving graph $M$ in an L-linkage, 
we observe that if the horizontal bar for edge $e$ is then outside of the range of the vertical stick for vertex $v$ 
--- where $(v,e)$ is a collision pair for $M$ --- then the collision which might have been caused by this collision
pair is perfectly avoided. However, those collision pairs of $M$ that do not fulfill this condition, still lead to collisions.
Inspired by this fact, we define 
the collision-freeness of an L-model as follows.
\begin{definition}\label{def:collision_free_L_model}
 Let $L=(M,h)$ be an L-model, where $M=(G,F)$ and $G=(V,E)$. It is {\bf collision-free} if and only if 
 $$h(e')\notin [\min_{v'\in e\in E}{h(e)},\max_{v'\in e\in E}{h(e)}]$$
 holds for any collision pair $(v',e')$ in $M$. If a collision pair fulfills this condition,
 we say that it is {\bf safe} (under function $h$). Hence, moving graph $M$ has a collision-free L-model if and 
 only if it has a height function $h$ such that all collision pairs of $M$ are safe under $h$.
\end{definition}

In this paper, we focus on the following problem: 
{\bf Given a finite moving graph, how to find a collision-free L-model for it.}
In the engineering context, this can be expressed as:
{\bf Given a finite moving graph, how to construct a collision-free
L-linkage for it.}

\section{Collecting collision pairs}\label{sec:collision_detection}

\begin{figure}[H]
\centering
\includegraphics[width=0.4\linewidth]{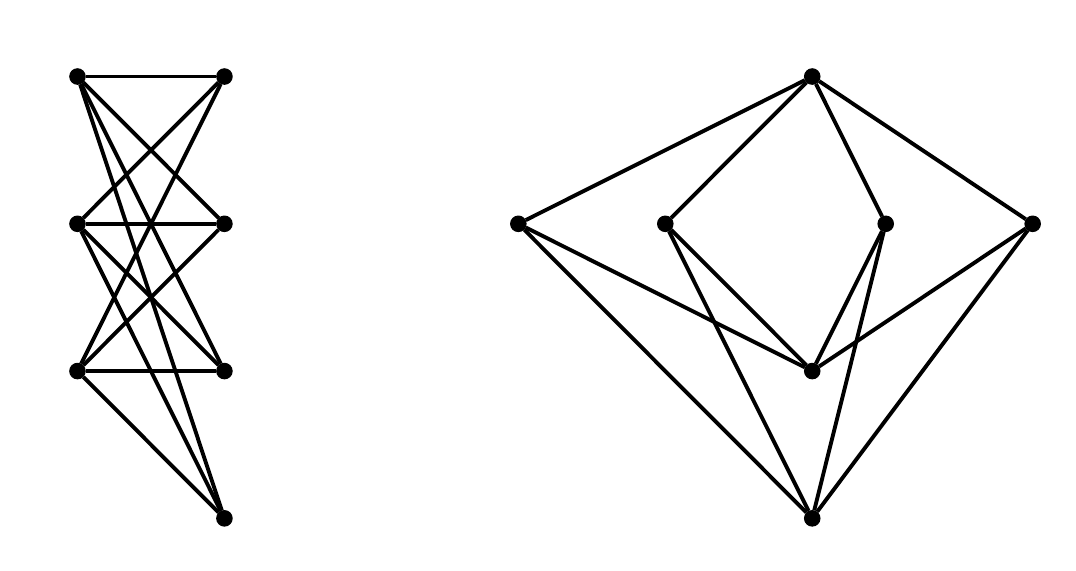}
\caption{The left sub-figure is a complete bipartite graph with 
seven vertices, while two independent sets contain three and four 
vertices, separately. The right sub-figure is the same graph,
while vertices of the two independent sets are transformed such that they lie on two orthogonal lines.
This graph is denoted by $K_{3,4}$ (or $K_{4,3}$).}
\label{fig:bipartite}
\end{figure}

A graph is {\em complete} if there is an edge
between any pair of vertices. An {\em independent set} of graph $G=(V,E)$ is a subset $V'$ of the vertex 
set $V$ such that there is no edge between vertices in $V'$. 
We say that graph $G=(V,E)$ is {\em bipartite}
if and only if there exists a bipartition of the vertex set $V=V_1\cup V_2$
such that both $V_1$ and $V_2$ are independent sets. And we call $V_1$ and $V_2$ {\em the two corresponding independent sets}
of the bipartite graph~$G$. 
We say that bipartite graph $G$ with two independent sets $V_1$ and $V_2$ is {\em complete}
if there is an edge connecting any two vertices $v,u$ for 
$v\in V_1$, $u\in V_2$. We denote by $K_{m,n}$ the complete bipartite graph 
with $m$, $n$ the cardinalities of its two independent sets respectively.
See Figure \ref{fig:bipartite} for an example of complete bipartite graph. 

Let $M=(G,F)$ be a moving graph with $G=(V,E)$. We detect collisions of this moving graph just by 
collecting its collision pairs. We check whether $(v,e)$ is a collision pair for all
$v\in V$ and all $e\in E$. Collection of all collision pairs reflects the collision information
of $M$. We implemented this program in Mathematica \cite{Mathematica}. Now let us see an example for a clearer idea.

\begin{example}\label{eg:collision_pairs}
 Let $G=(V,E)$ be the complete bipartite graph $K_{4,3}$ with vertex set $V=\{1,2,3,4,5,6,7\}$
 and two independent sets $V_1=\{1,2,3,4\}$, $V_2=\{5,6,7\}$.
 The moving graph $M=(G,F)$ belongs to the family of Dixon-1 moving graphs \cite{Dixon},
 where $F$ consists of the following functions:
 \begin{flalign*}
 &f_1(t) = (\sin{t},0), &&\\\nonumber
 &f_{2}(t) = (\sqrt{1+\sin^2{t}},0), &&\\\nonumber
 &f_{3}(t) = (-\sqrt{2+\sin^2{t}},0), &&\\\nonumber
 &f_{4}(t) = (\sqrt{3+\sin^2{t}},0), &&\\\nonumber
 &f_{5}(t) = (0,\cos{t}), &&\\\nonumber
 &f_{6}(t) = (0, \sqrt{1+\cos^2{t}}),&&\\\nonumber
 &f_{7}(t) = (0, -\sqrt{2+\cos^2{t}}). &&
 \end{flalign*}
 We proceed according to the definition of collision pairs. For each vertex $u$ and edge $\{v,w\}$, we check 
 if $f_u(t)$ lying on the line segment defined by 
 $f_v(t)$ and $f_w(t)$ has a solution in $\mathbb{R}$: if it does, then $(u,\{v,w\})$ is a collision pair of $M$.
 By going through all pairs of vertex $v\in V$ and edge $e\in E$ (such that $v\notin e$), we get the collision pairs in $M$ as follows:
$$(1,\{5,2\}), (1,\{5,3\}), (1,\{5,4\}), (2,\{5,4\}), (5,\{6,1\}), (5,\{7,1\}).$$ 

 \end{example}
 
With this we conclude the process of collecting collision pairs, which is the first step 
on our way of trying to find a collision-free L-model for the given moving graph.

\section{The partition condition}\label{sec:sufficient_condition}
In this section, we explain a sufficient condition for a moving graph to have a collision-free
L-model, which we call {\em the partition condition}. In order to introduce this condition,
some preparations are needed. Let us gradually approach it.

After collecting the collision pairs of a given moving graph, 
we want to express this information in a nicer way --- by a ``collision 
graph''.
\begin{definition}
Let $M=(G,F)$ be a moving graph, where $G=(V,E)$. The {\bf collision graph} $C=(V_C,E_C)$ of $M$ is 
a directed graph such that $V_C=E$ and $\overrightarrow{e_i,e_j}\in E_C$ if and only if at least one 
of the vertices of $e_i$, denoted by $v$, collides with
edge $e_j$ in $M$, i.e., $(v, e_j)$ is a collision pair of $M$, where $v\in e_i$ is one of the two incident
vertices of $e_i$.
\end{definition}

An {\em induced subgraph} $H=(V_1,E_1)$ of a given graph $G=(V,E)$ is another graph such that
$V_1\subset V$ and $E_1$ equals to the restriction of $E$ on $V_1$. That is to say,
edge set $E_1$ is formed of all edges of $G$ between any pair of vertices in $V_1$. 
Then $H$ is denoted by $G[V_1]$.
An induced collision graph
 is an induced subgraph of some given collision graph.

\begin{definition}[induced collision graphs]
 Let $M=(G,F)$ be a moving graph, where $G=(V,E)$. Let $S\subset E$, then the {\bf collision graph of $M$ induced by $S$}, denoted by $C[S]$, is 
 the subgraph of $C$ induced by $S$, where $C$ is the 
 collision graph of $M$.
\end{definition}

\begin{figure}[H]
\centering
\includegraphics[width=0.5\linewidth]{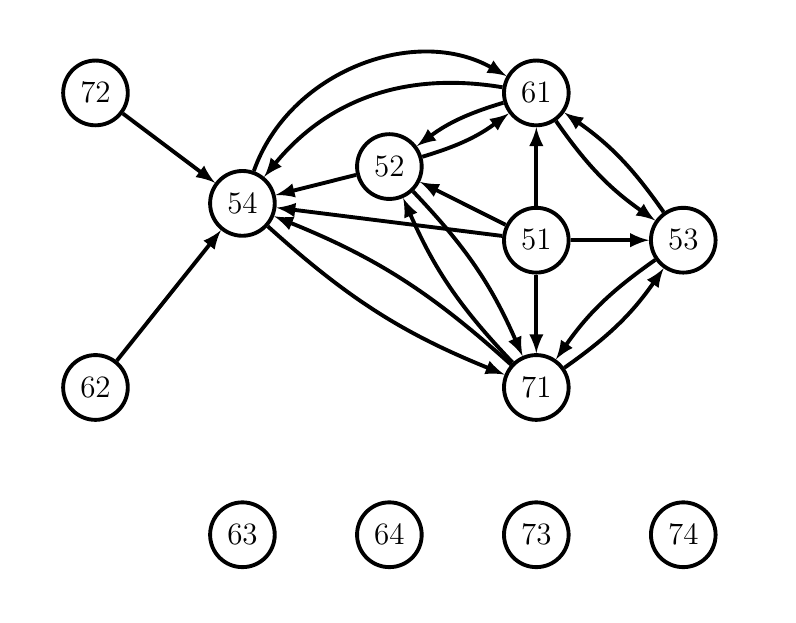}
\caption{This is the collision graph $C$ of the moving graph $M$ given in Example \ref{eg:collision_pairs}.
A directed edge $\{i,j\} \rightarrow \{k,l\}$ means either $(i,\{k,l\})$
or $(j,\{k,l\})$ is a collision pair of~$M$. Note that in the figure we write $ij$ for 
edge $\{i,j\}$, in order to have a more neat picture.}
\label{fig:collision_graph}
\end{figure}

  We continue with Example \ref{eg:collision_pairs}. From the collision information collected in Example~\ref{eg:collision_pairs}, we construct the 
  collision graph $C$, as shown in Figure \ref{fig:collision_graph}.
  
\begin{definition}
Let $C=(V_C,E_C)$ be a directed graph. If there exists a bipartition of $V_C$ into 
$V_1$ and $V_2$ such that $C[V_1]$ and $C[V_2]$ are both acyclic, then 
we say that graph $C$ fulfills {\bf the partition condition}.
\end{definition}

We say that a moving graph {\itshape fulfills the partition condition} if and only if its collision graph
fulfills the partition condition.
Now we are prepared for the main theorem.
\begin{theorem}\label{thm:partition_condition}
Let $M=(G,F)$ be a finite moving graph, where $G=(V,E)$. If $M$ fulfills the 
partition condition, then it has a 
collision-free L-model.
\end{theorem}

In order to prove Theorem \ref{thm:partition_condition}, we introduce an algorithm (see Algorithm \ref{alg:height_construction}) constructing 
the height function for the vertices of a given acyclic directed graph.
We need some more preparations before the proof. First we introduce an order
for the vertices of a directed graph.

\begin{definition}
Let $C$ be an acyclic directed graph. The {\bf height order} on the vertex set of $C$ --- denoted 
by ``$<$'' ---  is defined as:
$i<j$ if and only if there is a (directed) path from $i$ to $j$.
\end{definition}

\begin{remark}
A strict partial order is a binary relation that is irreflexive, transitive and asymmetric.
 The height order defined above is a strict partial order. 
\end{remark}

\begin{proposition}\label{prop:height_order}
 For a finite acyclic directed graph $C=(V,E)$, there exist(s) at least one minimal element 
 in $V$ with respect to the height order.
\end{proposition}

\begin{proof}
 If there was no minimal element in $V$ under the height order, then there would be an infinite chain $v_1>v_2>...$ in $V$.
 Since there is no cycle in graph $C$, 
 the elements
 in this chain are pairwise distinct. This contradicts the finiteness of the cardinality of~$V$.
\end{proof}

\begin{algorithm}[H]

\caption{Constructing the height function for elements in $V_C$.}\label{alg:height_construction}
\thispagestyle{empty}
\SetKwInOut{Input}{input}
\SetKwInOut{Output}{output}

\Input{a finite acyclic directed graph $C=(V_C,E_C)$; the height parameter 
    $(k_0,i)$ }
\Output{the height function $h:V_C\to Z$}
$k\gets k_0$ \;

\While{$C$ is not a null graph }
   {$S:=$ collection of all minimal vertices in the vertex set of graph $C$ 
   under height order\;
    $C:=C[V_C\setminus S]$, where $C[V_C\setminus S]$ is the subgraph of $C$ 
     induced by $V_C\setminus S$\;
     \While{$S\neq \emptyset$}
            {Pick one element $r$ in $S$, set $h(r):= k$\;
             $S:=S\setminus \{r\}$\;
             $k:=k+i$\;}}
\Return $h$.

\end{algorithm}
\begin{remark}
 Null graph is a graph with no vertices or edges.
 \end{remark}
 \begin{remark}\label{rem:output_not_unique}
 Note that the output of Algorithm \ref{alg:height_construction} is not unique, due to the fact that there can be 
 more than one minimal vertices under the height order in the considered directed graph, in each outer loop.
 This then leads to non-deterministic height assignments in the corresponding inner loops.
\end{remark}

In Algorithm \ref{alg:height_construction}, we add the height parameter as part of input data. The usage of this will be seen later on --- we want to
apply this algorithm to two directed graphs, but with two different height parameters, in order to 
construct a collision-free height function for the given moving graph. 
Note that Algorithm \ref{alg:height_construction} completes the height order which is a partial order, into a total order.
Now we can prove the termination of Algorithm \ref{alg:height_construction}.
\begin{proof}[Termination of Algorithm \ref{alg:height_construction}]
 If $C$ is not a null graph, by Proposition \ref{prop:height_order} there must exist at least 
 one minimal vertex under the height order. Hence $S$ is
 non-empty as long as $C$ is not null. Then the step ``$C:=C[V_C\setminus S]$'' strictly reduces the number of vertices of graph $C$, when $C$ is not a 
 null graph. 
 Since the given graph is finite, this algorithm
 terminates.
\end{proof}

The following proposition spells out the essence for the proof of Theorem \ref{thm:partition_condition}.
\begin{proposition}\label{prop:model_construction}
 Let $M=(G,F)$ be a moving graph that fulfills the partition condition, where $G=(V,E)$. 
 Let $C=(V_C,E_C)$ be the collision graph of $M$ and let $C_L=(V_L,E_L)$ and $C_U=(V_U,E_U)$ be the two corresponding acyclic induced collision 
 graphs (of $C$). After applying Algorithm \ref{alg:height_construction} with height parameter 
 $(1,1)$ to $C_U$, and with height parameter 
 $(0,-1)$ to $C_L$, we get a collision-free L-model of $M$.
\end{proposition}
\begin{remark}
 The above mentioned process assigns to edges (of $M$) in $V_U$ positive-integer height values, with the principle that $h(e_1)<h(e_2)$
if $\overrightarrow{e_1,e_2}\in E_U$. Meanwhile, it assigns to edges (of $M$) in $V_L$ non-positive-integer height values, with 
the principle that $h(e_1)>h(e_2)$ if $\overrightarrow{e_1,e_2}\in E_L$. 
\end{remark}

\begin{proof}[Proof of Proposition \ref{prop:model_construction}.]
Recall that vertices of the collision graph $C$ are exactly the edges of $M$.
In the sequel, whenever we talk about edges or vertices, we mean those with respect to the moving graph 
$M=(G,F)$, or the graph $G$ --- so as to make things less confusing.

Let $(v,e)$ be any collision pair of $M$. Then, $\overrightarrow{e_1,e}\in E_C$ if and only if $v\in e_1$, where $e_1$ is any edge of $M$. 
Without loss of generality, assume 
that $e\in V_U$ --- the case for $e\in V_L$ can be argued analogously. Considering the distribution of vertex~$v$ as one of the two 
incident vertices of any edge $e_1$ in the graphs $C_U$ and $C_L$, given that $(v,e)$ is a collision pair.
There are in total three cases:
\begin{itemize}
 \item First, $v$ shows up only in some edge(s) in graph $C_U$. In this case, Algorithm \ref{alg:height_construction} guarantees that the height 
 value of $e_1$ is strictly less than that of edge $e$, for all
 $\overrightarrow{e_1,e}\in E_U$. Hence, height value of edge $e$ is outside of the range of the height of edges (of $M$)
 containing $v$. That is to say, $(v,e)$ is safe.
 
 \item Second, $v$ shows up only in some edge(s) in graph $C_L$. In this case, Algorithm~\ref{alg:height_construction} guarantees 
 that the height values of all edges containing vertex $v$ are non-positive. Meanwhile, the height value of edge $e$ is positive. 
 Hence, we obtain that $h(e)\notin [\min_{v\in e'\in E}{h(e')}, \max_{v\in e'\in E}{h(e')}]$, $(v,e)$ is safe.
 
 \item Third, $v$ shows up in some edge(s) in graph $C_U$ and in some edge(s) in graph $C_L$ as well. From the above analysis, we know that
 edges containing $v$ that are in $V_L$ have non-positive height values and edges containing $v$ that are in $V_U$ has strictly less
 height values than edge $e$. Therefore, we have $h(e)\notin [\min_{v\in e'\in E}{h(e')}, \max_{v\in e'\in E}{h(e')}]$, $(v,e)$ is safe.
\end{itemize}
Hence we obtain that all collision pairs in $M$ are safe. By Definition 
\ref{def:collision_free_L_model}, indeed
we obtained a collision-free L-model of $M$ after the stated process.
\end{proof}

\begin{proof}[Proof of Theorem \ref{thm:partition_condition}.]
Straightforward, from Proposition \ref{prop:model_construction}.
\end{proof}
\begin{remark}
 From Proposition \ref{prop:model_construction}, we see that the height parameter in Algorithm \ref{alg:height_construction} is useful, 
 since we need to construct the edges of $M$ in $C_U$ from height $1$ upwards, while those in $C_L$
 are constructed from height $0$ downwards.
\end{remark}

For a clearer understanding towards Algorithm \ref{alg:height_construction}, we apply it to the moving graph $M$ given in 
Example \ref{eg:collision_pairs}, so as to get a collision-free L-model.

\begin{example}\label{eg:apply_height_algorithm}
We partition the edges of $M$ into two parts: $$E_U=\{\{5,1\},\{5,2\},\{5,3\},\{5,4\}\},$$
$$E_L=\{\{6,1\},\{6,2\},\{6,3\},\{6,4\},\{7,1\},\{7,2\}
,\{7,3\},\{7,4\}\},$$ then from these
two parts we construct the induced collision graphs $C_U$ and $C_L$, individually --- see Figure \ref{fig:induced_collision_graphs}.

  \begin{figure}[H]
\centering
\includegraphics[width=0.55\linewidth]{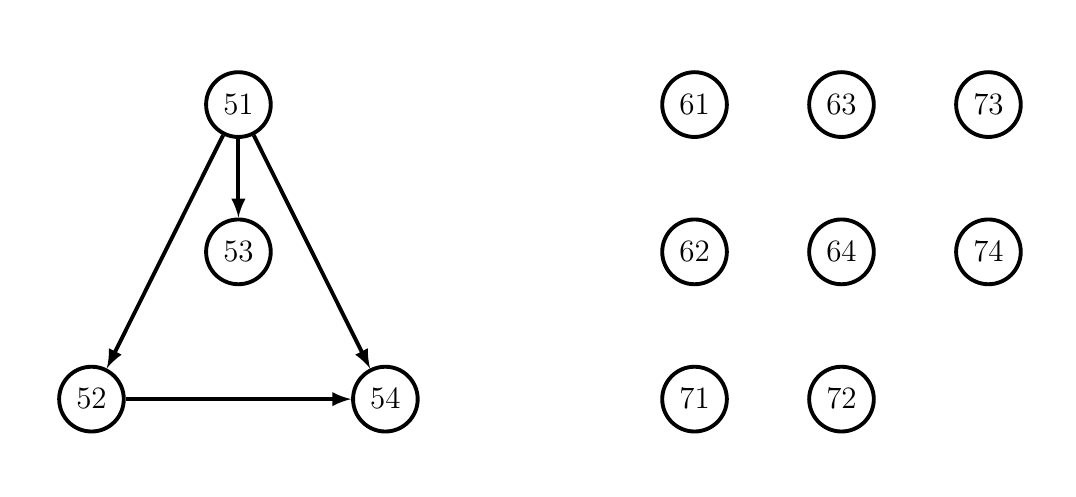}

\caption{The left sub-figure is the collision graph of $M$ described in Example \ref{eg:collision_pairs} 
induced by edge set $\{\{5,1\},\{5,2\},\{5,3\},\{5,4\}\}$ and the right sub-figure is the collision graph of $M$ induced by edge set
$\{\{6,1\},\{6,2\},\{6,3\},\{6,4\},\{7,1\},\{7,2\}
,\{7,3\},\{7,4\}\}$. Note that in the figure we write $ij$ for 
edge $\{i,j\}$, for convenience.}
\label{fig:induced_collision_graphs}
\end{figure}

  After applying Algorithm \ref{alg:height_construction}, according to Proposition \ref{prop:model_construction}, we get a height arrangement for some L-model as: 
  $h(\{5,4\})=4$, $h(\{5,3\})=3$, $h(\{5,2\})=2$, $h(\{5,1\})=1$, $h(\{6,1\})=0$, $h(\{6,2\})=-1$, $h(\{6,3\})=-2$, $h(\{6,4\})=-3$, 
  $h(\{7,1\})=-4$, $h(\{7,2\})=-5$, $h(\{7,3\})=-6$, $h(\{7,4\})=-7$.

  One can check that it indeed provides a collision-free L-model of $M$.
\end{example}

Note that the converse statement of Theorem \ref{thm:partition_condition} does not hold.
That is to say, there exists a moving graph $M=(G,F)$, where $G=(V,E)$, 
 which has a collision-free L-model but 
 we cannot partition $E$ into two parts $E_1$ and $E_2$ such that the induced
 collision graphs of $M$ (respectively by $E_1$ and $E_2$) are both acyclic.
We show it with the following example. 

\begin{example}\label{eg:converse_is_not_true}
  
 A moving graph fulfilling the above described conditions could be the one denoted by $S_2$ moving graph in
 \cite[Corollary 4.8]{moving_graph}.
 The formal definition of this moving graph (denoted by $M$ in our paper) is as follows.
 Let $M=(G,F)$, where $G=(V,E)$, $V=\{1,2,3,4,5,6,7,8\}$ and
 $$E=\{\{1,2\},\{1,4\},\{1,5\},\{3,2\},\{3,4\},\{3,5\},$$
 $$\{8,2\},\{8,4\},\{8,5\},\{1,7\},\{6,7\},\{6,5\},\{6,4\}\}.$$ 
 And $F$ consists of the following functions, where $a=1$, $b=\frac{11}{5}$, $c=\frac{3}{2}$:
\begin{flalign*}
 &s_1=\sqrt{b^2-a^2\cdot \sin^2{t}}\,, &&\\\nonumber
 &s_2=\sqrt{c^2-a^2\cdot \cos^2{t}}\,, &&\\\nonumber
 &f_1(t)=(-a\cdot \cos{t}-s_1, -a\cdot \sin{t}-s_2), &&\\\nonumber
 &f_2(t)=(a\cdot \cos{t}-s_1, -a\cdot \sin{t}+s_2), &&\\\nonumber
 &f_3(t)=(a\cdot\cos{t}+s_1, a\cdot\sin{t}+s_2), &&\\\nonumber
 &f_4(t)=(-a\cdot\cos{t}+s_1, -a\cdot\sin{t}+s_2), &&\\\nonumber
 &f_5(t)=(-a\cdot\cos{t}+s_1, a\cdot\sin{t}-s_2), &&\\\nonumber
 &f_6(t)=(-3\cdot a\cdot\cos{t}+s_1, -a\cdot\sin{t}-s_2), &&\\\nonumber
 &f_7(t)=(-3\cdot a\cdot\cos{t}-s_1, -a\cdot\sin{t}-3\cdot s_2), &&\\\nonumber
 &f_8(t)=(-a\cdot\cos{t}-s_1, a\cdot\sin{t}+s_2). &&
 \end{flalign*}
 After the collision-detection process, we get the following collision pairs:
 $$(2, \{8,4\}),\; (2, \{8,5\}),\;(3, \{1,4\}),\;(3, \{6,4\}),\; (3, \{8,4\}),\;(4, \{3,2\}),\; (4, \{3,5\}),\;$$
 $$(6, \{1,5\}),\;(6, \{3,5\}),\; (6, \{8,5\}),\;(5, \{6,4\}),\; (5, \{6,7\}),\;  (8, \{1,2\}),\; (8, \{3,2\}).$$

 One can check that the following height value assignments provide a 
 collision-free L-model for $M$:
 $$h(\{6,7\})=0,\; h(\{1,4\})=1,\; h(\{6,4\})=2,\; h(\{8,4\})=3,\; h(\{3,4\})=4,$$
 $$h(\{6,5\})=5,\; h(\{8,5\})=6,\; h(\{1,5\})=7,\; h(\{1,7\})=8,\; h(\{3,5\})=9,$$
 $$h(\{8,2\})=10,\; h(\{3,2\})=11,\; h(\{1,2\})=12.$$ 
 
 However, from the collision information of $M$, we know that $\overrightarrow{e_1,e_2}$,
 $\overrightarrow{e_2,e_1}$, $\overrightarrow{e_1,e_3}$, $\overrightarrow{e_3,e_1}$,
  $\overrightarrow{e_2,e_3}$ and $\overrightarrow{e_3,e_2}$ are all in $E_C$, where $e_1=\{3,2\}, e_2=\{6,4\}, e_3=\{8,5\}$ 
 and $E_C$ denotes the edge set of the collision graph $C=(V_C,E_C)$ of $M$. Hence, no matter how we try to 
 bipartition $V_C$ (i.e., $E$), we will get two among the edges $e_1,e_2,e_3$ in one 
 group, say $e_1$ and $e_2$. Then the edges $\overrightarrow{e_1,e_2}$ and 
 $\overrightarrow{e_2,e_1}$ already form a cycle in the corresponding induced collision graph.
 
 Therefore, we cannot partition $V_C$ into two parts such 
 that the induced collision graphs respectively of the two parts  
 are both acyclic.

\end{example}

To sum up this section, we propose the following steps, on how to decide whether a moving graph has a collision-free L-model.
\begin{enumerate}
 \item Collect the collision information of the given moving graph $M$, namely the collision pairs in $M$.
 \item Construct the collision graph of this moving graph, denote it by $C$.
 \item Decide whether $C$ fulfills the partition condition, or equivalently,
 whether $M$ fulfills the partition condition. If yes, with our
 algorithm we can get a collision-free L-model for $M$.
 
\end{enumerate}
\begin{remark}
 If the answer to Step 3. is ``no'', we can still apply a 
 brute-force algorithm which in worst case needs the factorial of $|E|$ many steps. 
 That is to say, try all possible (relative) height arrangements
 and check for each of them if it leads to a collision-free 
 L-model.   
 
\end{remark}

 But how to decide whether $C$ fulfills the partition condition in Step 3.? 
 To make the theory complete, now we introduce a method
 which may improve the efficiency in some cases, compare to the brute-force method.
 We need the definition of ``cycled subgraph'' in order to explain this
 method. Note that the word ``cycled'' here just refers to the situation
when there exist two edges between two vertices $v_1$ and $v_2$, namely 
$\overrightarrow{v_1,v_2}$ and $\overrightarrow{v_2,v_1}$, in the directed graph. In this case,
these two edges form a cycle in graph $C$. Therefore, these two vertices must stay in different 
induced collision graphs --- if $C$ fulfills the partition condition.
We use this graph to depict such information.

\begin{definition}
 Let $C=(V,E)$ be a directed graph. Collect all vertices $v,u\in V$
 such that $\overrightarrow{v,u}$ and $\overrightarrow{u,v}$ 
 are both edges of $C$ in set $S$. Then consider the graph $C[S]$. First, delete all the single edges.
 Then, view all edges between any two vertices as one single non-directed edge.
 We call the so-obtained graph the {\bf cycled subgraph} of graph $C$.
\end{definition}

\begin{example}
 Let $C$ denote the graph depicted in Figure \ref{fig:collision_graph}.
 The cycled subgraph $C_1$ of $C$ is shown in 
 Figure \ref{fig:cycled_subgraph}.

  \begin{figure}[H]
\centering
\includegraphics[width=0.4\linewidth]{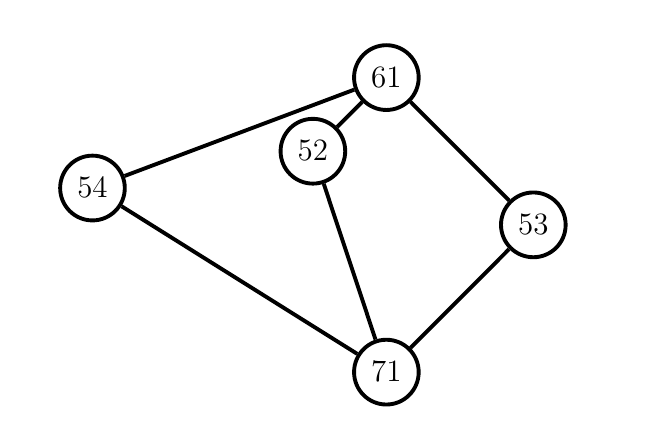}

\caption{This is the cycled subgraph $C_1$ of graph $C$ which
is shown in Figure \ref{fig:collision_graph}.}
\label{fig:cycled_subgraph}
\end{figure}

\end{example}

 It is not hard 
to see that the following proposition holds since any edge in the cycled subgraph
indicates a cycle in the given directed graph.

\begin{proposition}\label{prop:cycled_subgraph_bipartite}
 If the cycled subgraph $S=(V_S,E_S)$ of a directed graph $C=(V,E)$ is
 not bipartite, then the graph $C$ does not 
 satisfy the partition condition. Any 
 bipartition of $V=V_L\cup V_U$ such that both $C[V_L]$
 and $C[V_U]$ are acyclic induces a bipartition of $V_S$ which 
 gives the two independent sets of graph $S$. 
 \end{proposition}
 \begin{proof}
  Suppose that graph $C$ satisfies the partition condition and the two corresponding 
  induced subgraphs are $C_L=(V_L,E_L)$ and $C_U=(V_U,E_U)$. If there exist $v_1,v_2\in V_S\cap V_L$
  such that $\{v_1,v_2\}\in E_S$,
  by the definition of cycled subgraph we know that the two edges $\overrightarrow{v_1,v_2}$
  and $\overrightarrow{v_2,v_1}$ are both in $E_L$. This makes $C_L$ cyclic, which is a contradiction.
  Hence, $\{v_1,v_2\}\notin E_S$. Analogously, we have that if $v_1,v_2\in V_S\cap V_U$,
  then $\{v_1,v_2\}\notin E_S$. Hence, $S$ is bipartite, the two independent sets of $S$ are exactly $V_S\cap V_L$
  and $V_S\cap V_U$ --- note that $V_S\subset V_L\cup V_U=V$.
 \end{proof}

 \begin{remark}
 We see that a bipartition of the edges of the graph in Figure \ref{fig:collision_graph} such that the two obtained induced collision graphs 
 are both acyclic, say the one stated
 in Example~\ref{eg:apply_height_algorithm}, includes a bipartition
 on the edge set $$\{\{5,2\},\{5,3\},\{5,4\},\{6,1\},\{7,1\}\}$$ of graph $C_1$ drawn in Figure~\ref{fig:cycled_subgraph}
  as follows: 
 $$ \{\{5,2\},\{5,3\},\{5,4\}\}\cup \{\{6,1\},\{7,1\}\},$$ which also indicates
 a bipartite structure of graph $C_1$ --- there are no edges in between vertices 
 $\{\{5,2\},\{5,3\},\{5,4\}\}$ or $\{\{6,1\},\{7,1\}\}$.
 \end{remark}
 \begin{remark}
  Consider the $S_2$ moving graph. From Example~\ref{eg:converse_is_not_true}, we
  know that the cycled subgraph of $S_2$ moving graph contains a triangle. Therefore, it cannot be bipartite, 
  hence the $S_2$ moving graph does not fulfill the partition condition.
 \end{remark}

 We hope that this, as a pre-step of the algorithm for deciding whether a given directed graph
 fulfills the partition condition, can save labor for us in some situation. Now we give an algorithm 
 (see Algorithm \ref{alg:partition_condition}) for deciding whether a given 
 directed graph fulfills the partition condition.

 \begin{algorithm}

\caption{Deciding whether a directed graph $C$ fulfills the partition condition}\label{alg:partition_condition}
\thispagestyle{empty}
\SetKwInOut{Input}{input}
\SetKwInOut{Output}{output}

\Input{a finite directed graph $C=(V_C,E_C)$ }
\Output{If $C$ fulfills the partition condition, output the corresponding bipartition of $V_C$; otherwise, output ``No''.}
$S=(V_S,E_S)\gets$ the cycled subgraph of $C$ \;
 \If{$S$ is not bipartite}
     {\Return ``No''}
 \Else 
  {Find all bipartitions of $V_S$ such that both parts are independent sets of $S$\;
    Denote them by the sequence 
    $B_1,...,B_m$, where $B_i:=\{p_i,q_i\}$ for $1\leq i\leq m$\;
    List all bipartitions of $V_C\setminus V_S$ as ordered pairs (hence each bipartition is considered twice) 
    and denote the pairs by the sequence $B'_1,...,B'_k$, where $B'_i:=(p'_i,q'_i)$
    for $1\leq i\leq k$\;
   \For{$i$ from $1$ to $m$}
     {\For{$j$ from $1$ to $k$}
          {Check whether $(p_i\cup p'_j,q_i\cup q'_j)$ 
           is a partition of $V_C$ such that the induced 
           subgraph of both parts are acyclic\;
           \If{yes}
               {\Return $(p_i\cup p'_j,q_i\cup q'_j)$}}}
   \Return ``No'' }
  
\end{algorithm}
  
\begin{remark}
 In Algorithm \ref{alg:partition_condition}, we consider each bipartition of $V_C\setminus V_S$ twice, i.e., consider them as 
 ordered pairs, while the bipartitions of $V_S$ are viewed as cardinality-two sets, considered only once. This is to make sure that we 
 indeed go through all combinations of the bipartitions of $V_S$ and those of $V_C\setminus V_S$, once.
\end{remark}

\begin{proof}[Termination and correctness of Algorithm \ref{alg:partition_condition}]
 Termination is obvious since the input is a finite graph. Correctness
 follows from Proposition \ref{prop:cycled_subgraph_bipartite}.
\end{proof}

 When the collision information for some moving graph contains two vertices colliding with each other, 
 then we already know that the range of these two vertices should not intersect in the L-model (so as to avoid collisions). 
 This information can also save us 
 some work when we try to find a collision-free height arrangement, but we will not go too much into details on it. 
 
  Note that our approach assumes that a moving graph, with the parameterization of the positions of the vertices in 
 a variable $t$ is given. Computing such a parameterization for a given graph with fixed edge lengths
 is a preliminary step, which in practice 
 should be addressed before seeking a collision-free L-model for the moving graph. 
 In this preliminary step, the corresponding configuration space is 
 an algebraic curve in many cases. In general, parameterization of algebraic curves is a topic on its own 
 and it is not the focus of this paper. There are many references available on this topic.
 Some details on algorithmic parameterization of rational curves can be checked in 
 \cite[Chapter~4]{rational_curves}, for instance.  Examples of some computations on the traces of vertices for visualizing the 
 moving graphs (via animations) are available in \cite{computation}. 
 On the parameterization step for the two families of moving graphs that are discussed in this paper, we provide 
 the following ideas:
 Formulas for Dixon-1 moving graphs can be obtained and
 verified by Pythagorean Theorem. Formulas for Dixon-2 moving graphs can be deduced by its
 defining equation system and its symmetry property.

 With this we complete the description of our main results. In the coming 
 sections, we apply the theory to two classical families of moving graphs:
 Dixon-1 and Dixon-2 moving graphs.

\section{Dixon-1 moving graphs}

In this section we focus on Dixon-1 moving graph family. The
Dixon-1 construction can be applied
to an arbitrary bipartite graph, giving a mobile realization; the so-obtained moving graph would be a {\em Dixon-1 moving graph}.
In the sequel, we explain this construction in natural language.
Let $G=(V,E)$ be an arbitrary 
bipartite graph and let $V_1$, $V_2$ be its two independent sets. 
Then, place the vertices in $V_1$ to the $x$-axis and those in $V_2$ to the
$y$-axis, where by $x$-axis and $y$-axis we just refer to two orthogonal axes in the real plane.
Note that the positions of vertices are not fixed, the only requirement is that the
vertices in $V_1$ should stay on the $x$-axis, while those in $V_2$ should stay on the $y$-axis.
One can check that the dimension of the configuration space is one for complete bipartite graphs. That is to say,
we have one-degree of freedom to choose the position of vertices. It is a mobile realization of complete
bipartite graphs, and also of all bipartite graphs in the cases of which then may have higher dimensional configuration spaces. 
See the right sub-figure of Figure \ref{fig:bipartite} for one such placement of the 
vertices.

A more specific definition (with a fixed scale) for Dixon-1 moving graphs would be the following.
Let $m$ and $n$ be positive integers and $K_{m,n}$ the complete 
undirected bipartite graph with 
$m+n$ vertices. Fix real numbers $0<a_1<...<a_{m-1}$ and $0<b_1<...<b_{n-1}$,
$a_0=b_0=0$. The two independent sets of vertices are $P:=\{p_0,...,p_{m-1}\}$ and $Q:=\{q_0,...,q_{n-1}\}$. 
We place $p_i$ ($0\leq i\leq m-1$) to $x$-axis and $q_j$ ($0\leq j\leq n-1$) to $y$-axis. 
The positions of the vertices of graph $K_{m,n}$ parameterized in the variable $t$ are as follows:
\begin{flalign*}
 &f_{p_0}(t) = (\sin{t},\; 0), &&\\\nonumber
 &f_{p_i}(t) = (\pm\sqrt{a_i+\sin^2{t}},\; 0),\;  i> 0, &&\\\nonumber
 &f_{q_0}(t) = (0,\; \cos{t}), &&\\\nonumber
 &f_{q_j}(t) = (0,\; \pm\sqrt{b_j+\cos^2{t}}),\; j> 0. &&
 \end{flalign*}
 Note that we can choose the $x$-coordinate for $f_{p_i}(t)$ ($i>0$) to be either positive 
 or negative; the same rule applies to the $y$-coordinate for $f_{q_j}(t)$ ($j> 0$).
We denote the first coordinate (x-coordinate) of the position of vertex $v$ by $f_v(t)|_x$ and the second coordinate 
(y-coordinate) of it by $f_v(t)|_y$. We observe that the edge lengths are fixed. 
We compute the lengths for edges $\{p_i, q_j\}$ ($0\leq i \leq m-1$, $0\leq j \leq n-1$) as follows:
 $$\|f_{p_i}(t)-f_{q_j}(t)\|=\sqrt{a_i+\sin^2{t}+b_j+\cos^2{t}}=\sqrt{a_i+b_j+1}.$$
 
 The next theorem tells us the collision information of an arbitrary Dixon-1 moving 
 graph.
\begin{theorem}\label{thm:dixon1}
 The collision pairs of any Dixon-1 moving graph $D$
 (defined as above) are as follows:
 \begin{enumerate}
  \item $(p_0,\{p_i,q_0\})\in CP_D$, $i> 0$.

  \item $(p_i, \{p_k, q_0\})\in CP_D$, for $k>i>0$ and $f_{p_k}(t)|_x\cdot f_{p_i}(t)|_x>0$.

  \item $(q_0, \{p_0, q_i\})\in CP_D$, $i> 0$.

  \item $(q_i, \{p_0, q_k\})\in CP_D$, for $k>i>0$ and $f_{q_k}(t)|_y\cdot f_{q_i}(t)|_y>0$.

 \end{enumerate}

\end{theorem}

\begin{proof}
\begin{enumerate}
 \item In the case when $f_{p_i}(t)=(\sqrt{a_i+\sin^2{t}},\; 0)$: Solving the equation 
 $$\|f_{q_0}(t)-f_{p_0}(t)\|+\|f_{p_0}(t)-f_{p_i}(t)\|=\|f_{q_0}(t)-f_{p_i}(t)\|,$$
i.e., 
$$1+\sqrt{a_i+\sin^2{t}}-\sin{t}=\sqrt{a_i+1},$$
we obtain that $t=\frac{\pi}{2}+2\cdot k\cdot \pi$ 
($k\in \mathbb{Z}$) are in the solution set. 

In the case when $f_{p_i}(t)=(-\sqrt{a_i+\sin^2{t}},\; 0)$: Solving the equation 
 $$\|f_{q_0}(t)-f_{p_0}(t)\|+\|f_{p_0}(t)-f_{p_i}(t)\|=\|f_{q_0}(t)-f_{p_i}(t)\|,$$
i.e., 
$$1+\sin{t}+\sqrt{a_i+\sin^2{t}}=\sqrt{a_i+1},$$
we obtain that $t=-\frac{\pi}{2}+2\cdot k\cdot \pi$ 
($k\in \mathbb{Z}$) are in the solution set. 
By definition of collision pairs, $p_0$ does not collide with
edge $\{p_0,q_0\}$.
Hence, $(p_0,\{p_i, q_0\})\in CP_D$, $i> 0$.

Are these already all the 
collision pairs containing vertex $p_0$? Suppose that $p_0$
collides with edge $\{p_k, q_l\}$, $k,l>0$, then there exists $t'\in \mathbb{R}$
such that $f_{p_0}(t')$ lies on the line segment defined by the 
points $f_{p_k}(t')$ and $f_{q_l}(t')$. Since 
$f_{p_0}(t')|_y=f_{p_k}(t')|_y=0$, we have 
$$f_{q_l}(t')|_y= \sqrt{b_l+\cos^2{t}}=0$$
or 
$$f_{q_l}(t')|_y= -\sqrt{b_l+\cos^2{t}}=0.$$
However, $b_l>0$ for $l>0$, so 
neither of the above equations hold. We reach a contradiction.
Therefore, $(p_0,\{p_i, q_0\})$ ($i>0$) are already all the collision pairs 
of $D$ that have $p_0$ as the vertex.

\item In the case when both $f_{p_i}(t)|_x$ and $f_{p_k}(t)|_x$ are positive:
Solving the equation $\|f_{q_0}(t)-f_{p_i}(t)\|+
\|f_{p_i}(t)-f_{p_k}(t)\|=\|f_{q_0}(t)-f_{p_k}(t)\|$, i.e.,
$$\sqrt{a_i+1}+\sqrt{a_k+\sin^2{t}}-\sqrt{a_i+\sin^2{t}}=\sqrt{a_k+1},$$
we see that $t=\frac{\pi}{2}+k_1\cdot \pi$, $k_1\in \mathbb{Z}$ are 
in the solution set. 

When $f_{p_i}(t)|_x<0$ and $f_{p_k}(t)|_x<0$, the equation is the same,
hence we obtain the same set of solutions. 
But are these already all the 
collision pairs containing vertex $p_i$ ($i>0$)? Suppose that $p_i$
collides with edge $\{p_k, q_l\}$ for some $k<i$. Then there exists 
$t'\in \mathbb{R}$ such that the following equation holds:
$$\|f_{q_l}(t')-f_{p_i}(t')\|+\|f_{p_i}(t')-f_{p_k}(t')\|=\|f_{q_l}(t')-f_{p_k}(t')\|.$$
We obtain that $\|f_{q_l}(t')-f_{p_i}(t')\|\leq \|f_{q_l}(t')-f_{p_k}(t')\|$, i.e.,
$\sqrt{a_i+b_l+1}\leq \sqrt{a_k+b_l+1}$. However, $0\leq a_k<a_i$ since $k<i$. We gained a contradiction.
So, $(p_i, \{p_k,q_l\})\notin CP_D$ when $k<i$. 

Suppose that $p_i$ collides with edge $\{p_k,q_l\}$ for $k>i>0$.
There exists $t'\in\mathbb{R}$ such that $f_{p_i}(t')$ lies on the line segment defined by 
$f_{q_l}(t')$ and $f_{p_k}(t')$. Since the $y$-coordinates of both $f_{p_k}(t')$ and $f_{p_i}(t')$
are zero, we obtain that $f_{q_l}(t')|_y=0$. This leads to the restriction that $l=0$, since
 $f_{q_l}(t')|_y\neq 0$ when $l> 0$. 
 
 Suppose that $p_i$ collides with edge $\{p_k, q_0\}$ for $k>i>0$ and the corresponding $t'\in \mathbb{R}$
 makes $f_{p_k}(t')|_x\cdot f_{p_i}(t')|_x<0$. Suppose that $f_{p_k}(t')|_x>0$ and $f_{p_i}(t')|_x<0$.
 Since both $f_{q_0}(t')|_x$ and $f_{p_k}(t')|_x$ are non-negative and $f_{p_i}(t')$ lies on the line 
 segment defined by these two points, we obtain that $f_{p_i}(t')|_x\geq 0$. This contradicts the assumption $f_{p_i}(t')|_x<0$.
 In the case when $f_{p_k}(t')|_x<0$ and $f_{p_i}(t')|_x>0$, the argument is analogous.
 
 To sum up, $(p_i, \{p_k, q_0 \})$ ($k>i$, $f_{p_k}(t)|_x\cdot f_{p_i}(t)|_x>0$) are already all the collision pairs 
 containing vertex $p_i$ ($i>0$).

 \item Because of the symmetry of vertices in $Q$ and vertices in $P$ in the sense of exchanging
 $x$-axis with $y$-axis, the remaining proofs can be done analogously. The proof for item 3. is analogous to 
 that of item 1.
 \item The proof for item 4. is analogous to that of item 2.\phantom\qedhere \qed
 \end{enumerate}
\end{proof}

After some analysis on the collision information, we find out that we can apply 
the method introduced in Section \ref{sec:sufficient_condition} and then 
obtain a collision-free L-model for any Dixon-1 graph! Before we can talk about the main
theorem of this section, we need to introduce a lemma, before which however, we need 
two definitions first.

  Let $D_1=(K_{m,n}, F_1)$ be a Dixon-1 moving graph with two independent sets 
  $P:=\{p_0,\ldots, p_{m-1}\}$ and
 $Q:=\{q_0,\ldots, q_{n-1}\}$ such that
 all $x$-coordinates of the positions of vertices in $\underline{P}:=P\setminus\{p_0\}$ have the same sign and all
 $y$-coordinates of the positions of vertices in $\underline{Q}:=Q\setminus\{q_0\}$ have the same sign. Such Dixon-1
 moving graphs are called {\bf crowded}. 
 Let $D_2=(K_{m,n}, F_2)$ be any Dixon-1 moving graph with the same underlying graph $K_{m,n}$
 as $D_1$. Let $D_2$ be so that the only difference between $D_1$ and $D_2$ is on the difference 
 between $F_1$ and $F_2$ --- on the choices of the signs for $x$-coordinates of vertices 
 in $\underline{P}$ and $y$-coordinates of vertices in $\underline{Q}$. Then, we say that $D_1$ is the 
 {\bf crowded extreme} of $D_2$.

\begin{lemma}\label{lem:crowded}
Let $D_1$ be the crowded extreme of $D_2$, where $D_1$, $D_2$ are two Dixon-1 moving graphs defined as above.
Then, we have $CP_{D_2}\subset CP_{D_1}$.
\end{lemma}
\begin{proof}
 By Theorem \ref{thm:dixon1} item 1., we have that
 $(p_0,\{p_i, q_0\})\in CP_{D_2}$ and $(p_0,\{p_i,q_0\})\in CP_{D_1}$ for $i>0$.
 By Theorem \ref{thm:dixon1} item 2., we have that 
 $(p_i, \{p_k, q_0\})\in CP_{D_2}$, for $k>i>0$ and $f_{p_k}(t)|_x\cdot f_{p_i}(t)|_x>0$.
 However, since all $x$-coordinates of the positions of vertices in $\underline{P}$ have the same sign in $D_1$,
 all these collision pairs in $CP_{D_2}$ are in $CP_{D_1}$ as well.
 The analyses for the collision pairs in $CP_{D_2}$ described in item 3. and item 4. 
 being in $CP_{D_1}$ are analogous.
 Therefore, $CP_{D_2}\subset CP_{D_1}$.
\end{proof}

\begin{corollary}\label{cor:crowded_extreme}
 Let $D$ be an arbitrary Dixon-1 moving graph and let $h$ be a height function
 of $D$ such that $(D,h)$ is a collision-free L-model. Let $D'$
 be any Dixon-1 moving graph of which the crowded extreme is $D$, then 
 $(D',h)$ is a collision-free L-model of $D'$.
\end{corollary}
\begin{proof}
 From Lemma \ref{lem:crowded}, we know that all collision pairs of $D'$ are
 collision pairs of $D$. Since $h$ guarantees that all collision pairs 
 in $D$ are safe,
 $h$ naturally guarantees that all collision pairs in $D'$ are safe. 
 Therefore, $(D',h)$ is a 
 collision-free L-model of $D'$.
 \end{proof}

\begin{theorem}\label{thm:dixon1_has_model}
 Every Dixon-1 moving graph has a collision-free L-model.
\end{theorem}
\begin{proof}
 By Corollary \ref{cor:crowded_extreme}, it suffices to show the existence of a collision-free L-model
 for any crowded Dixon-1 moving graph. Let $D=(K_{m,n}, F)$, $K_{m,n}=(V,E)$ be any crowded Dixon-1 moving graph
 with independent sets $P:=\{p_0,\ldots, p_{m-1}\}$ and $Q:=\{q_0,\ldots, q_{n-1}\}$. Recall that the 
 edge set of $D$ is 
 $$E=\{\{p_i, q_j\}\mid 0\leq i\leq m-1, 0\leq j\leq n-1 \}.$$
 Let $C=(V_C,E_C)$ be the collision graph of $D$, recall that $V_C=E$.
 In the sequel, we give a bipartition of the set $E$, into $E_1$ and $E_2$,
 such that the two corresponding induced collision graphs $C[E_1]$ and
 $C[E_2]$ of $D$ are both acyclic. We divide the edges in $E$ into two groups as follows:
 \begin{itemize}
\item $E_1=\{\{p_i,q_0\}\mid 0\leq i\leq m-1\}$;
 \item   $E_2=E\setminus E_1$.
 \end{itemize}
By Theorem \ref{thm:dixon1} item 1., 
there is a directed edge from $\{p_0,q_0\}$ to each of the other elements in $E_1$ 
in graph $C[E_1]$. By Theorem \ref{thm:dixon1} item 2., there is a directed
edge from $\{p_i,q_0\}$ to $\{p_k,q_0\}$ in graph $C[E_1]$ for all pairs $(i,k)$ such that $i<k$. Theorem \ref{thm:dixon1}
tells us that these are all the edges
in graph $C[E_1]$. We easily see that the graph is acyclic.

By Theorem \ref{thm:dixon1} item 4., there is a directed edge from 
$\{p_j, q_i\}$ to $\{p_0, q_k\}$ for $k>i>0$ and $j\geq 0$ in graph $C[E_2]$. These
are all the edges of graph $C[E_2]$. Since these edges in $C[E_2]$ are directed from some 
edge of $D$ with strictly smaller subscript on the $q$-vertex to some edge
with strictly bigger subscript on the $q$-vertex, there is no cycle in
graph $C[E_2]$.

By Theorem \ref{thm:partition_condition}, $D$ has a collision-free L-model. 
As Proposition \ref{prop:model_construction} stated, we apply Algorithm \ref{alg:height_construction} with height parameters
$(1,1)$ to $C[E_1]$ and $(0,-1)$ to $C[E_2]$, respectively. Then we obtain 
$$h(\{p_i, q_0\})=i+1,$$
$$h(\{p_l, q_j\})=-(j-1)\cdot (m+1)-l.$$ 
This is the corresponding
height function for a collision-free L-model of $D$; moreover, it induces
a collision-free L-model for any Dixon-1 moving graph which has $D$ as its 
crowded extreme as well, by Corollary \ref{cor:crowded_extreme}.
\end{proof}

\begin{example}
 Example \ref{eg:apply_height_algorithm} uses exactly the same grouping (of edges) as stated in the proof 
 of Theorem
 \ref{thm:dixon1_has_model}. And the height function obtained coincides with the result as if we compute using the above formulae. 
 Let us again greet this running example.

 Let $M=(K_{3,4},F)$, where $F$ is the set of functions given in Example \ref{eg:collision_pairs}. 
 Two independent sets of $K_{3,4}$ are $P=\{p_0=1,p_1=2,p_2=3,p_3=4\}$ 
 and $Q=\{q_0=5,q_1=6,q_2=7\}$. Then it is not hard to check that 
 the grouping of edges and
 the height function obtained in Example \ref{eg:apply_height_algorithm} indeed coincide with those provided
 in the proof of Theorem \ref{thm:dixon1_has_model}. 
 
 However, note that there are also other possible output 
 height functions of Algorithm~\ref{alg:height_construction} which provide collision-free L-models of $M$,
 by Remark \ref{rem:output_not_unique}.
\end{example}

\section{Dixon-2 moving graphs}

In this section we discuss another class of moving graphs --- Dixon-2 moving graph family. 
 The Dixon-1 construction
applies to any bipartite graph, however, the Dixon-2 construction 
only applies to the complete bipartite graph $K_{4,4}$ and its subgraphs.

\begin{figure}[H]
\centering
\includegraphics[width=0.35\linewidth]{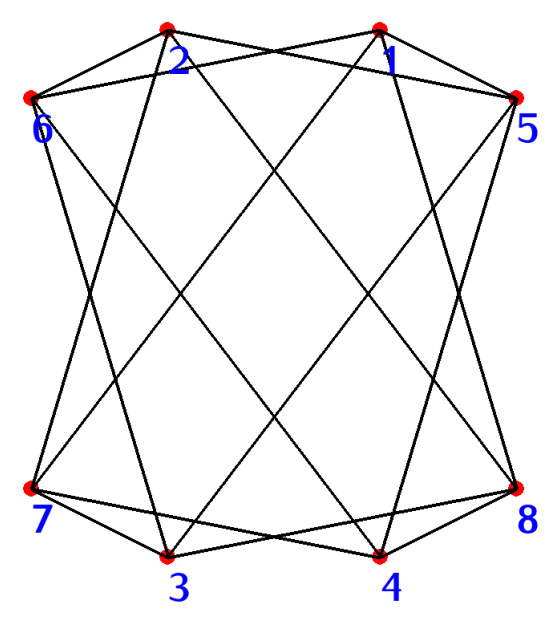}
\caption{This is a static Dixon-2 moving graph.}
\label{fig:dixon2}
\end{figure}

To say it in natural language:
Fix vertices $1$ and $5$ in the first quadrant of the coordinate system. 
Let vertex $2$ be symmetric to vertex $1$ with respect to $y$-axis, let vertex $4$ be symmetric 
to vertex $1$ w.r.t. $x$-axis, and let vertex $3$ be symmetric to vertex $1$ w.r.t. the origin.
Let vertex $6$ be symmetric to vertex $5$ w.r.t. $y$-axis, let vertex $8$ be symmetric 
to vertex $5$ w.r.t. $x$-axis, and let vertex $7$ be symmetric to vertex $5$ w.r.t. the origin. 
 Note that here $x$-axis refers to the horizontal axis and $y$-axis refers to the vertical axis, as in the convention.
 Let Group One consist of vertices $1,2,3,4$ and Group Two constitute of vertices $5,6,7,8$.
Then we add an edge between any
vertex in group one and any vertex in group two, forming a bipartite graph $K_{4,4}$. This is the Dixon-2 construction. 
To apply it to any subgraph,
simply take the substructure of this construction. If we apply this construction to $K_{4,4}$, we then obtain a 
{\em Dixon-2 moving graph}. See Figure \ref{fig:dixon2} for a visualization. 

Some other literature may include also those moving graphs obtained via
applying the construction to proper subgraphs of $K_{4,4}$ in the Dixon-2 moving graph families; our notion differs from them slightly.
Now we give a more formal definition for Dixon-2 moving graphs as follows.
Let the vertex set of $K_{4,4}$ be
$\{1,2,3,4,5,6,7,8\}$ and the two independent sets be $\{1,2,3,4\}$ and $\{5,6,7,8\}$.
Let $b,c,d\in \mathbb{R}^+$ be such that $d>a$ and $b>a$ hold. The parameterizations 
of the positions of the vertices are as follows:
 \begin{flalign*}
 &s_1 = \sqrt{b^2-a^2\cdot \sin^2{t}},\;\;
 s_2 = \sqrt{d^2-a^2\cdot\cos^2{t}}\;, &&\\\nonumber
 &f_{1}(t) = \left(\frac{a\cdot \cos{t}+s_1}{2},\; \frac{a\cdot \sin{t}+s_2}{2}\right), &&\\\nonumber
 &f_{2}(t) = \left(-\frac{a\cdot \cos{t}+s_1}{2},\; \frac{a\cdot \sin{t}+s_2}{2}\right), &&\\\nonumber
 &f_{3}(t) = \left(-\frac{a\cdot \cos{t}+s_1}{2},\; -\frac{a\cdot \sin{t}+s_2}{2}\right), &&\\\nonumber
 &f_{4}(t) = \left(\frac{a\cdot \cos{t}+s_1}{2},\; -\frac{a\cdot \sin{t}+s_2}{2}\right), &&\\\nonumber
 &f_{5}(t) = \left(\frac{-a\cdot \cos{t}+s_1}{2},\; \frac{-a\cdot \sin{t}+s_2}{2}\right), &&\\\nonumber
 &f_{6}(t) = \left(-\frac{-a\cdot \cos{t}+s_1}{2},\; \frac{-a\cdot \sin{t}+s_2}{2}\right), &&\\\nonumber
 &f_{7}(t) = \left(-\frac{-a\cdot \cos{t}+s_1}{2},\; -\frac{-a\cdot \sin{t}+s_2}{2}\right), &&\\\nonumber
 &f_{8}(t) = \left(\frac{-a\cdot \cos{t}+s_1}{2},\; -\frac{-a\cdot \sin{t}+s_2}{2}\right). &&
 \end{flalign*}

We denote the coordinates of vertex $1$ by $(x_1,y_1)$ and those of vertex $5$ by $(x_5, y_5)$.
Then from the above parameterizations one can check that $x_1,y_1,x_5,y_5 > 0$ holds.
We can acquire the coordinates of the other six vertices individually, by alternating the sign of one or both of the coordinates 
of either vertex $1$ or vertex $5$. 
Let $c\in \mathbb{R}^+$ such that $c^2=b^2+d^2-a^2$ holds. Since $b,d>a$, we know that $c>b>a$ and $c>d>a$. 
Then we have the following four equations, where we require $x_1,y_1,x_5,y_5 > 0$ additionally. 

\begin{align*}
 (x_1-x_5)^2+(y_1-y_5)^2 = a^2 \\
 (x_1+x_5)^2+(y_1-y_5)^2 = b^2 \\
 (x_1+x_5)^2+(y_1+y_5)^2 = c^2 \\
 (x_1-x_5)^2+(y_1+y_5)^2 = d^2 
\end{align*}
Since $a^2+c^2=b^2+d^2$, one can verify that the third equation is redundant. Actually, the above equation system has one-dimensional
 solution set, hence the realization is mobile. We see that the edge lengths are fixed as follows, where $l_{u,v}:=\|f_u(t)-f_v(t)\|$ denotes the length of edge $\{u,v\}$:
 $$l_{1,5}=l_{2,6}=l_{3,7}=l_{4,8}=a,\; l_{1,6}=l_{2,5}=l_{3,8}=l_{4,7}=b,$$ 
 $$l_{1,7}=l_{2,8}=l_{3,5}=l_{4,6}=c,\; l_{1,8}=l_{2,7}=l_{3,6}=l_{4,5}=d.$$
 
 Since we require $x_1,y_1,x_5,y_5>0$, vertices $1$ and $5$ stay in the first quadrant; 
 vertices $2$ and $6$ stay in the second quadrant; vertices $3$ and $7$ stay in the third quadrant and
 vertices $4$ and $8$ stay in the fourth quadrant. We call the vertices that are in the same quadrant
 {\bf friend vertex} of each other.
 
 Above, we define the Dixon-2 moving graphs in two different ways: one is via 
 parameterizing the positions of the vertices using explicit functions on $t$, the other is by giving the equations that the coordinates
 of the vertices need to fulfill using implicit functions.
 Actually, these two ways are equivalent. To consider the collision-detecting equations, 
we can use both systems --- the implicit one, or the explicit one. 
In the upcoming theorem, we analyze the collision information for an arbitrary Dixon-2 moving graph.
\begin{theorem}
The collision pairs of any Dixon-2 moving graph $D$
 (defined as above) are as follows --- written in $8$ groups each of which illustrates the collision pairs of one vertex:
 \begin{enumerate}
  \item $(1,\{5,2\})$, $(1,\{5,4\})$, $(1,\{5,3\})$;
  \item $(5,\{1,6\})$, $(5, \{1,8\})$, $(5,\{1,7\})$;
  \item $(2,\{6,1\})$, $(2,\{6,3\})$, $(2,\{6,4\})$;
  \item $(6,\{2,5\})$, $(6,\{2,7\})$, $(6,\{2,8\})$;
  \item $(3,\{7,1\})$, $(3,\{7,4\})$, $(3,\{7,2\})$;
  \item $(7,\{3,5\})$, $(7,\{3,8\})$, $(7,\{3,6\})$;
  \item $(4,\{8,1\})$, $(4,\{8,2\})$, $(4,\{8,3\})$;
  \item $(8,\{4,5\})$, $(8,\{4,6\})$, $(8,\{4,7\})$.
 \end{enumerate}
 To sum up, elements in $CP_D$ are $(i,\{k,l\})$, where $1\leq i\leq 8$, $k$ is the friend vertex of $i$,
 and $l$ is in the same independent set with $i$.
 
\end{theorem}
\begin{proof}
Because of the symmetry property of Dixon-2 construction, it suffices to verify the collision situation for vertex $1$. 
We will show that vertex $1$ collides with all edges containing vertex $5$ --- its friend vertex --- except for the edge $\{1,5\}$, of course;
and we will show that these are all the collision pairs containing vertex $1$.

Solving the equation 
$$\|f_1(t)-f_5(t)\|+\|f_1(t)-f_2(t)\|=\|f_5(t)-f_2(t)\|$$ 
which is equivalent to 
$$\|(x_1,y_1)-(x_5,y_5)\|+\|(x_1,y_1)-(-x_1,y_1)\|=\|(x_5,y_5)-(-x_1,y_1)\|,$$
i.e., $a+2\cdot x_1=b$, we obtain that
$$x_1 = \frac{b-a}{2},\; x_5 = \frac{a+b}{2},\; 
 y_1 = y_5=\frac{\sqrt{d^2-a^2}}{2}$$ which is equivalent to $t=\pi+2\cdot k\cdot \pi$, $k\in \mathbb{Z}$.
Hence $(1,\{5,2\})\in CP_D$. Solving the equation 
$$\|f_1(t)-f_5(t)\|+\|f_1(t)-f_4(t)\|=\|f_5(t)-f_4(t)\|$$ 
which is equivalent to 
$$\|(x_1,y_1)-(x_5,y_5)\|+\|(x_1,y_1)-(x_1,-y_1)\|=\|(x_5,y_5)-(x_1,-y_1)\|,$$
i.e., $a+2\cdot y_1=d$, we obtain that
$$y_1 = \frac{d-a}{2},\; y_5 = \frac{a+d}{2},\; 
 x_1 = x_5=\frac{\sqrt{b^2-a^2}}{2}$$ which is equivalent to $t=-\frac{\pi}{2}+2\cdot k\cdot \pi$, $k\in \mathbb{Z}$.
Hence $(1,\{5,4\})\in CP_D$. 
Solving the equation 
$$\|f_1(t)-f_5(t)\|+\|f_1(t)-f_3(t)\|=\|f_5(t)-f_3(t)\|$$
which is equivalent to 
$$\|(x_1,y_1)-(x_5,y_5)\|+\|(x_1,y_1)-(-x_1,-y_1)\|=\|(x_5,y_5)-(-x_1,-y_1)\|,$$
i.e., $a+\sqrt{(2\cdot x_1)^2+(2\cdot y_1)^2}=c$, we get
 $$x_1=\frac{1}{2}\sqrt{\frac{(b^2-a^2)(c-a)}{c+a}},\; x_5=\frac{1}{2}\sqrt{\frac{(b^2-a^2)(c+a)}{c-a}},$$
 $$y_1=\frac{1}{2}\sqrt{\frac{(c^2-b^2)(c-a)}{c+a}},\; y_5=\frac{1}{2}\sqrt{\frac{(c^2-b^2)(c+a)}{c-a}},$$
which is a solution fulfilling the implicit defining equation system of $D$ and it is not hard to verify that
 $x_1,y_1,x_5,y_5>0$ holds. Hence $(1,\{5,3\})\in CP_D$.

Since vertex $1$ by definition stays in the first quadrant, it is impossible for it to collide with edges 
that do not intersect the first quadrant. Therefore, the only edges remain to be considered for
vertex $1$ are $\{6,4\}$ and $\{2,8\}$. Because of the symmetry property of Dixon-2 construction, it suffices
to consider whether $(1,\{6,4\})$ is a collision pair or not. Solving the equation 
$$\|f_1(t)-f_6(t)\|+\|f_1(t)-f_4(t)\|=\|f_6(t)-f_4(t)\|$$
which is equivalent to 
$$\|(x_1,y_1)-(-x_5,y_5)\|+\|(x_1,y_1)-(x_1,-y_1)\|=\|(-x_5,y_5)-(x_1,-y_1)\|,$$
i.e., $a+2\cdot y_1=c$, we get $y_1=\frac{c-b}{2}$. Substituting it back to the explicit defining equation 
system of $D$, we obtain that $\sin{t}=-\frac{b}{a}$, which implies $|\sin{t}|>1$. This contradicts the 
fact that $0\leq |\sin{t}|\leq 1$. 
Hence $(1,\{6,4\})\notin CP_D$. There are in total three collision pairs of $D$ that contains vertex $1$, namely
$(1,\{5,2\})$, $(1,\{5,4\})$, and $(1,\{5,3\})$.

The collision situation for vertex $5$ is analogous. 
By the symmetry property of Dixon-2 construction, collision situations for the remaining 
vertices can be obtained directly. 
 \end{proof}

In the sequel, we prove that Dixon-2 moving graphs have no
collision-free L-models. 

\begin{theorem}\label{thm:dixon2}
There are no collision-free L-models for Dixon-2 moving graphs.
\end{theorem}
\begin{proof}
We prove by contradiction.
Let $D$ be an arbitrary Dixon-2 moving graph and let $h$ be a height function of $D$ such that $(D,h)$ is collision-free. 
Then, let us analyze the relative heights for edges of this moving graph.

 Since vertex $5$ collides with all edges containing vertex $1$ (except for edge $\{1,5\}$), 
 there are no edges containing vertex $1$ lying in between edge $\{1,5\}$ and another edge containing 
 vertex $5$. Symmetrically, there are no edges 
 containing vertex $5$ lying in between edge $\{1,5\}$ and another edge containing 
 vertex $1$. To conclude, the range of vertex $1$ and that of vertex $5$ 
 intersect on edge $\{1,5\}$. W.l.o.g., assume that edges containing vertex $1$ (except for edge $\{1,5\}$) 
 have bigger height values than those containing vertex~$5$ (except for edge $\{1,5\}$).
 Then we obtain 
 $$h(\{1,x\})>h(\{1,y\})>h(\{1,z\})>h(\{1,5\})>h(\{u,5\})>h(\{v,5\})>h(\{w,5\}),$$ 
 where $\{x,y,z\}=\{6,7,8\}$ and $\{u,v,w\}=\{2,3,4\}$. 
 
 Now we claim that 
 $h(\{4,7\})>h(\{3,5\})$. Suppose this does not hold, i.e., $h(\{4,7\})<h(\{3,5\})$. Then 
 we obtain $$h(\{1,7\})>h(\{3,5\})>h(\{4,7\}).$$ 
 Then we see that edge $\{3,5\}$ is within the range 
 of vertex $7$, this contradicts the fact that the collision pair $(7,\{3,5\})$ being safe. 
 Hence we have $h(\{4,7\})>h(\{3,5\})$. Now we claim that 
 $h(\{3,8\})>h(\{4,5\})$. Suppose this does not hold, i.e., $h(\{3,8\})<h(\{4,5\})$. Then 
 we obtain $$h(\{1,8\})>h(\{4,5\})>h(\{3,8\}).$$ 
 Then we see that edge $\{4,5\}$ is within the range 
 of vertex $8$, this contradicts the fact that the collision pair $(8,\{4,5\})$ being safe. 
 Therefore we obtain $h(\{3,8\})>h(\{4,5\})$. 
 
 Thereafter, we try to analyze the relative heights of the edges $\{4,7\}$ and $\{3,8\}$. If $h(\{4,7\})>h(\{3,8\})$,
 we have $$h(\{4,7\})>h(\{3,8\})>h(\{4,5\}).$$ 
 Then we see that edge $\{3,8\}$ is within the range 
 of vertex $4$, this contradicts the fact that the collision pair $(4,\{3,8\})$ being safe.
 If $h(\{4,7\})<h(\{3,8\})$, then
 we have $$h(\{3,5\})<h(\{4,7\})<h(\{3,8\}).$$
 Then we see that edge $\{4,7\}$ is within the range 
 of vertex $3$, this contradicts the fact that the collision pair $(3,\{4,7\})$ being safe.
 Hence, no matter how we place edge $\{4,7\}$ and edge $\{3,8\}$ relatively, it always leads
 to some collision pair being unsafe under $h$, which contradicts the assumption of $(D,h)$
 being collision-free.
 
 Hence, there are no collision-free L-models for any Dixon-2 moving graph. 
\end{proof}

As side products, we obtain the following two relatively trivial corollaries.

\begin{corollary}
 Dixon-2 moving graphs do not fulfill the partition condition. 
 That is to say, there is no partition of the edges of any Dixon-2 moving graph into 
 two parts $E_L$, $E_U$, such that the induced collision graphs $C_L$(by $E_L$) 
 and $C_U$(by $E_U$) are both acyclic.
\end{corollary}
\begin{proof}
 If a Dixon-2 moving graph fulfilled the partition condition, then by Theorem~\ref{thm:partition_condition}, it would have a collision-free L-model,
 which contradicts Theorem \ref{thm:dixon2}. 
\end{proof}

Naturally, we define a {\bf moving super-graph of $M$} to be a moving graph that contains the moving graph $M$ 
as a substructure. Upon this concept, we have the following corollary.
\begin{corollary}
Let $M=(G,F)$ with underlying graph $G=(V,E)$ be a moving super-graph of a Dixon-2 moving graph $D=(G_1,F_1)$ with underlying graph 
$G_1=(V_1,E_1)$.
Then, the moving graph $M$ does not have any collision-free L-models.
\end{corollary}
\begin{proof}
 Suppose that $(M,h')$ is a collision-free L-model. Since $F$ equals $F_1$ when restricted to the vertex set $V_1$, 
 any collision pair $(v,e)$
 of $D$ is also a collision pair of $M$. Then all collision pairs 
 $(v,e)$ of $D$ are safe under $h'$, since $(M,h')$ is a collision-free L-model. 
 Now we define a function $h:E_1\to \mathbb{Z}$ to be the restriction 
 of $h':E\to \mathbb{Z}$ on $E_1$. Then for any $(v,e)\in CP_D$, $h(e)$ should be outside of the height range of 
 vertex~$v$ in $D$; otherwise $h'(e)$ would also be within the range of $v$ in $M$, since
 the range of $v$ in $M$ is equal or wider than that in $D$. 
 Then, all collision pairs of $D$ are safe under $h$; this contradicts Theorem \ref{thm:dixon2}.
\end{proof}

We wind up the section by confirming that,
in reality we cannot build any collision-free Lego models (as shown in Figure 
\ref{fig:photo}) or L-linkages for any Dixon-2 moving graph.

\section{Acknowledgement}
This research was funded by the Austrian Science Fund (FWF): Grant no.W1214-N15, project DK9;
it was also supported by the strategic program ``Innovatives OÖ 2020'' by the Upper
Austrian Government. 

I am genuinely thankful to Josef Schicho for proposing this interesting problem to me, 
the parameterization of Dixon-2 moving graphs,
suggestions on the previous-work part of this paper, comments on modifications, and offering me the Python code 
for the animation of Dixon-2 moving graphs (from which Figure \ref{fig:dixon2}
was screen-shot).
 I am sincerely grateful to Jan Legerský
for reminding me of the edge intersection case, providing me with
the parameterization of $S_2$ moving graphs, especially
for all those down-to-each-sentence detailed comments on the previous manuscript, which to a great extend improve the quality of 
this paper. 
I thank Ralf Hemmecke for helping with debugging the collision-detection program.

\end{document}